# Bars and spheroids in gravimetry problem

## Valery Sizikov[1] and Vadim Evseev[2]

[1]University ITMO, Kronverksky pr., 49, 197101 Saint-Petersburg, Russia
[2] Palacký University Olomouc, tř. 17 listopadu, 1192/12, 771 46 Olomouc, Czech Republic

Email: sizikov2000@mail.ru and vadevs020913@gmail.com

**Abstract**

The direct gravimetry problem is solved by dividing each deposit body into a set of vertical adjoining bars, whereas in the inverse problem, each deposit body is modelled by a homogeneous ellipsoid of revolution (spheroid). Well-known formulae for the *z*-component of gravitational intensity for a spheroid are transformed to a convenient form. Parameters of a spheroid are determined by minimizing the Tikhonov smoothing functional with constraints on the parameters, which makes the ill-posed inverse problem by unique and stable. The Bulakh algorithm for initial estimating the depth and mass of a deposit is modified. The proposed technique is illustrated by numerical model examples of deposits in the form of two and five bodies. The inverse gravimetry problem is interpreted as a gravitational tomography problem or, in other words, as 'introscopy' of Earth's crust and mantle.

## 1. Introduction

Modelling of deposits is one of principal approaches used in solving the direct and, especially, inverse gravimetry problems. The *direct problem* is the computation of the gravitational field produced by some modelled deposit on Earth's surface. The *inverse problem* is the determination of the deposit parameters from the modelled or measured field anomaly (e.g., the Bouguer anomaly) on Earth's surface. For the direct problem, the approximation of a deposit by several bodies of arbitrary shape is one of preferable approaches. In the inverse problem, one should use the bodies of approximately regular form, but which corresponds to the shape of the deposit.

### 1.1. Models of deposits

Some authors use the following simplified *models of deposits* [7, 10, 37, 47]: in the form of quadrangular truncated pyramids, prisms, cylinders, beams, polyhedrons, parallelepipeds, intersecting bars, etc. However, such figures have non-smooth surfaces and generate cumbersome (although not complicated) formulae (see, e.g., [47]). Let us also mention the plane-layered model [42]. In works [17–19], homogeneous (and inhomogeneous) *spheroids*, or ellipsoids of revolution, are used as the deposit models and in work [20] et al., spheroids are used for modelling the Earth figure. Such models are effectively applied, e.g., in astrophysics for constructing the galactic models [43, 57]. In this work, we continue to use spheroids for deposit modelling.

More complex (algebraic) models have also been developed, namely, 3D models and the inversion of gravity data using discretization grids with thousands and even millions of cells involving the solution of large matrix systems (e.g., [10, 15, 22, 29, 61], et al.). Such approaches are more general and accurate, but more complex and require more computational resources. In some cases, the use of simplified models (e.g., in the form of spheroids) is more obvious and convenient as well as can serve as a good initial approximation when constructing complex models.

Different kinds of *initial data* is used for calculating deposits, namely, the intensities of gravitational and magnetic fields [4, 7, 17–19, 34, 41], gravity gradient tensor components [37, 61], seismic data [4, 27], remote sensing from satellites [4, 32, 41], etc. In this work, only the vertical component $V_z$ of the gravitational field intensity on Earth's surface is used, which is, in principle, sufficient for calculating rather complex spheroidal deposit models.

In works [1, 27, 28, 52, 56] et al., various methods of determining the boundary shape $z(x)$ separating the two parts of Earth's crust (properly the crust and deposit) are described. The measured function is the anomaly of the gravitational force $\Delta g(x)$ (e.g., the Bouguer anomaly [41]). In addition, the lower boundary of the deposit $H = $ const [52, 56] or the coordinates of its centre [1, 28] as well as the density anomaly $\Delta\rho$ are set. In this case, a one-dimensional nonlinear integral equation is solved with respect to $z(x)$ (the upper boundary of the deposit [52, 56] or even its entire boundary [1, 28]). Furthermore, each cross-section of the deposit is modelled by an ellipse [28] or the boundary is arbitrary [1]. In [1, 28, 52, 56], the problem is solved as a set of one-dimensional problems (in a number of vertical cross-sections). It is important, that a number of parameters is a priori set, namely, the lower boundary of the body $H$, the coordinates of its center and the density difference $\Delta\rho$. As a result, such a technique for solving the inverse problem is simplified, although the solution is unique.

In a number of papers [2, 12, 33, 36, 39, 53] et al., to calculate the mass potential (vanishing outside the Earth), the authors use harmonic (as well as anharmonic and biharmonic) functions and different constraints for the density $\rho$ (biharmonic constraint, the boundedness of the density $0 \leq \rho \leq \rho_{\max}$, etc.), that can provide uniqueness of the inverse problem.

In this work, the inverse problem is solved not by setting the deposit parameters but they are determined at the expense of decremental constraints and a regularization.

In this work in the *direct problem* (the calculation of the gravity field anomaly induced by the deposit on Earth's surface), a deposit is approximated by several bodies each of which is approximated by a set of adjoining *vertical bars* [18, 19]. In the *inverse problem* (determination of the deposit parameters from the field anomaly measured at Earth's surface or modelled in the direct problem), each body of the deposit is modelled by a *biaxial ellipsoid*, or *ellipsoid of revolution*, or *spheroid* (convenient astrophysical term).

Spheroids (ellipsoids) are widely used in celestial mechanics [11, 51], astrophysics (galactic models) [43, 57] and geophysics [59]. However, formulae for the field induced by an ellipsoid, as they are given in most of the works, can be reduced to a more convenient form. Moreover, after the works of Yun'kov [59], spheroids have not been often used in geophysics [17, 19, 20, 25, 50].

It is assumed in [50] that the bodies of the ore type being sources of gravitational field have shapes close to spheroids. Spheroids are homogeneous, convex and star-shaped domains having the mean plane. Uniqueness theorems hold for such domains [36, 46]. In the inverse problem in [50], the shapes of the bodies are determined by minimizing the discrepancy functional using the Lagrange undetermined multipliers (the regularization parameters). However, the coordinates of the centres and densities of the bodies have to be set, which limits the applicability of the method.

It is assumed in [25] that local inclusions have the shapes of homogeneous bodies of revolution, in particular, that of a sphere, oblate or prolate spheroid. The solution of the inverse problem is sought in the form of a series in polynomials and associated Legendre functions as well as in the form of splines by minimizing the Tikhonov functional with the aid of the variation method.

For the direct problem, it is assumed in this work that geologic bodies have rather arbitrary shapes (though resembling spheroids). In the inverse problem, the bodies are modelled by spheroids. Furthermore, in this work, the coordinates of the centres of the spheroids, their semi-axes and densities as well as the quantity of the spheroids are included into a number of the unknowns (with constraints on the values).

*1.2. Comparison with tomography*

The inverse gravimetry problem is often solved as a set of two-dimensional problems.

First, the density anomaly distribution is determined for some of the vertical cross-sections, and, second, a three-dimensional (volumetric) image is composed. This procedure quite resembles approaches typical for different types of tomography [30, 31, 44, 45, 58], and first of all, those of the X-ray computerized tomography (XCT) and magnetic resonance tomography (MRT) [31, 44, 45, 58]. In this work, which uses spheroids (as in [17–19]), the three-dimensional problem is solved which is equivalent to the three-dimensional tomography [9]. It is therefore proposed to refer the inverse gravimetry problem to the class of tomography problems (which was already done in [17–19]) and to call the inverse gravimetry problem the *gravitational tomography problem*. Moreover, a special feature of the considered problem is that it enables to 'look' inside Earth by means of mathematical processing instead of drilling wells which reminds of the process of introscopy being typical also of tomography. Therefore, it is also proposed to call the inverse gravimetry problem the *Earth introscopy problem* [17–19]. The proposed interpretations suggest the possibility of application of the extensive developments in the field of computerized tomography, in particular, XCT and MRT to solve the inverse gravimetry problem.

## 2. Calculation of the direct problem using a set of vertical bars

Let us consider a geologic deposit in the form of several homogeneous bodies having rather arbitrary shapes. As an example, a model of a deposit in the form of two bodies is shown in figure 1. In this example, body 1 is conditionally associated with an ore body and body 2 with an intrusion.

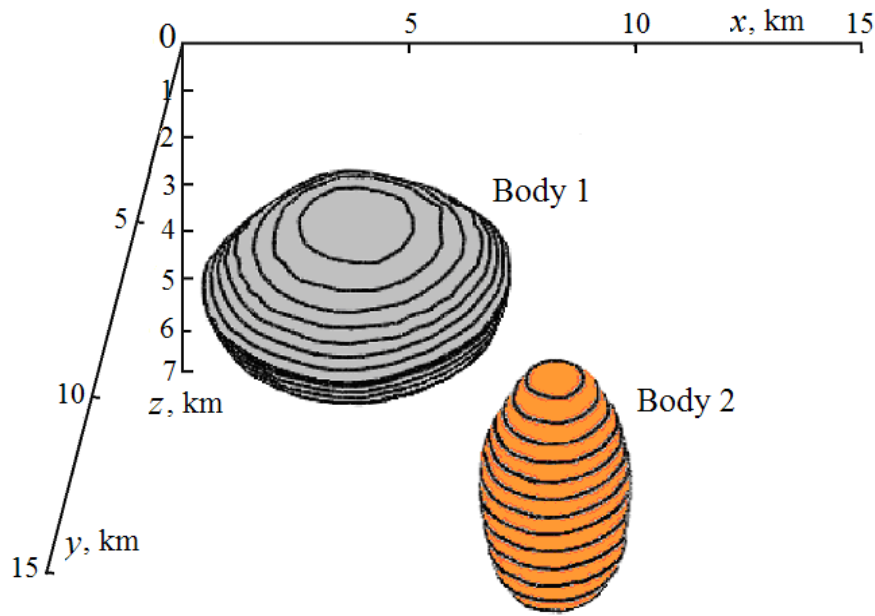

**Figure 1**. A model of deposit in the form of two bodies.

The contours of the $z$-sections for each body are shown in figure 2. The numbers on the contours are the $z$-coordinates (in km) of the $z$-sections. Furthermore, the continuous lines are the contours located above the conditional median section (the mean plane according to the geologic models of the Sretenskii class [5, 46]), and the short-dashed lines are the contours under that section.

**Definition 2.1** [18, 19]. A body is called *vertically star-shaped* if any vertical ray (straight line) intersects its boundary only twice.

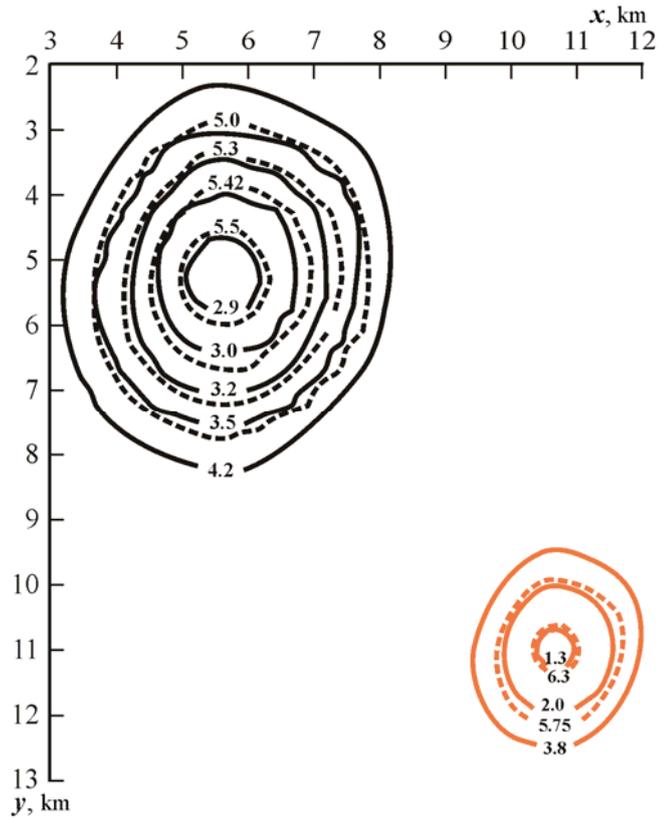

**Figure 2**. The contours of the $z$-sections of the modelled deposit bodies.

Let $\rho = \text{const}$ be the density of a body and let the body be vertically star-shaped. Let us approximate it by a set of elementary vertical *bars* with cross-sections $dx'dy'$ and with boundaries $z'_{\min} = z'_{\min}(x', y')$ and $z'_{\max} = z'_{\max}(x', y')$ (the roof and bottom according to the terminology of [5]) (figure 3). In the direct as well as inverse problems, only the $z$-components of the gravitational intensity induced by deposit bodies are considered.

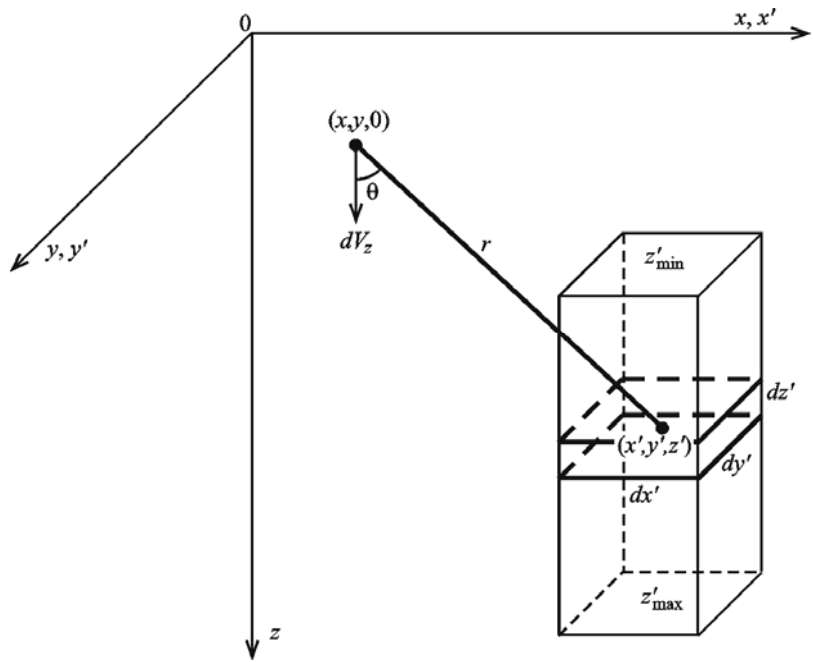

**Figure 3**. The elementary cell $dx'dy'dz'$ and elementary vertical bar of a deposit body.

**Lemma 2.1** (Gravitational field induced by an elementary bar). *The intensity $V_z$ induced by an elementary bar at a point $(x, y, 0)$ equals* [19]

$$dV_z(x, y, 0) = \gamma \rho \left[ \frac{1}{\sqrt{(x-x')^2 + (y-y')^2 + z'^2_{\min}}} - \frac{1}{\sqrt{(x-x')^2 + (y-y')^2 + z'^2_{\max}}} \right] dx' dy', \quad (2.1)$$

*where $\gamma$ is the gravitational constant and $x', y', z'$ are the coordinates of the elementary cell.*

**Proof.** The $z$-intensity induced by an elementary cell $dx'dy'dz'$ of a deposit body at a point $(x, y, 0)$ (figure 3) equals

$$dv_z = \gamma \frac{\rho \, dx' dy' dz'}{r^2} \cos\theta = \frac{\gamma \rho z'}{r^3} dx' dy' dz', \quad (2.2)$$

where $r = \sqrt{(x-x')^2 + (y-y')^2 + z'^2}$. Then, the $z$-intensity induced by an elementary bar at the point $(x, y, 0)$ is equal to the integral of the expression (2.2)

$$dV_z(x, y, 0) = \int_{z'_{\min}}^{z'_{\max}} dv_z = \gamma \rho \left[ \int_{z'_{\min}}^{z'_{\max}} \frac{z'}{r^3} dz' \right] dx' dy', \quad (2.3)$$

which yields the expression (2.1). □

**Corollary 2.1.** The integration of (2.1) over all elementary bars adjacent to each other from $z'_{\min} = z'_{\min}(x', y')$ to $z'_{\max} = z'_{\max}(x', y')$ yields the intensity $V_z(x, y, 0)$ induced by the entire body at the point $(x, y, 0)$. But if the condition for a body being vertically star-shaped is not satisfied, the 'voids' must be excluded from the integration regions $[z'_{\min}, z'_{\max}]$.

**Remark.** Strictly speaking, in the direct problem, a body is not modeled by a set of vertical bars. The bars are used only for the calculation of the field. Furthermore, a body can have an arbitrary shape.

This approach for modelling the direct problem is rather simple and effective, which is confirmed by the numerical examples (see the examples in [18, 19] and in this paper below).

## 3. Modelling the inverse problem using spheroids

**Definition 3.1.** An *ellipsoid* is a body bounded by the surface

$$\xi^2/a^2 + \eta^2/b^2 + \zeta^2/c^2 = 1,$$

where $a$, $b$ and $c$ are the semiaxes of the ellipsoid and the origin of the coordinate system is placed at the body centre.

In this work, biaxial ellipsoids, or ellipsoids of revolution around axis $z$, referred to as *spheroids* for short further in the text, are considered. In works [11, 20, 51, 59], formulae for the potential $V$ and intensity components $V_x$, $V_y$, $V_z$ of a spheroid have been deduced. However, formulae have not been reduced to a certain convenient form in the indicated works. This is attempted to be done further in this work for $V_z$.

**Remark.** Following usual practice [11], the term ellipsoid (spheroid) is used to refer to the body bounded by the surface as well as only the surface.

Consider an *oblate spheroid* for which $a = b > c$. In [11, 51, 59], a coordinate system $x$, $y$, $z$ was introduced with the origin at the centre of the spheroid and with the $z$-axis directed vertically upwards along the minor axis $c$ of the spheroid. For this case, formulae have been deduced for the potential $V(x, y, z)$ and intensity components $V_x(x, y, z)$, $V_y(x, y, z)$, $V_z(x, y, z)$ induced by the spheroid at a point $(x, y, z)$ outside the spheroid.

However, we will consider the case when the origin of the coordinate system $x, y, z$ is placed at some point on Earth's surface with the $z$-axis directed vertically downwards. Let $x_0, y_0, z_0$ be the coordinates of spheroid's centre and let the $z$-component $V_z$ of the field induced by the spheroid be measured (or computed) at a point $(x, y, 0)$ (figures 3 and 4).

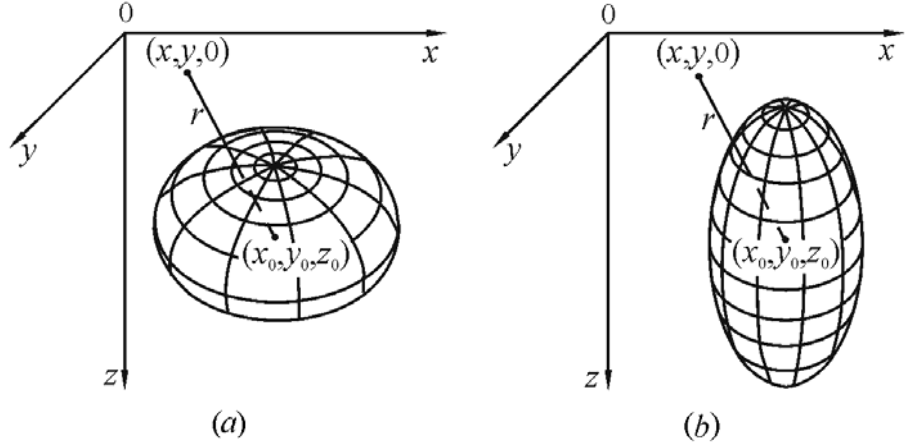

**Figure 4.** Oblate (*a*) and prolate (*b*) spheroids.

### 3.1. Formulae for z-intensity induced by a spheroid

**Lemma 3.1.** *For an oblate (along z) spheroid (figure 4a), the formula for the z-intensity induced at a point* $(x, y, 0)$ *outside the spheroid takes the form:*

$$V_z(x, y, 0) = 4\pi\gamma\rho \frac{\varepsilon}{e^3}(p - \operatorname{arctg} p)\, z_0, \qquad (3.1)$$

*where*

$$\varepsilon = c/a < 1, \quad e = \sqrt{1 - \varepsilon^2} > 0, \quad p = q/\sqrt{\tau}, \quad q = ea/r,$$
$$r = \sqrt{(x - x_0)^2 + (y - y_0)^2 + z_0^2}, \quad \tau = \left[1 - q^2 + \sqrt{(1 - q^2)^2 + 4q^2 z_0^2/r^2}\right]/2. \qquad (3.2)$$

**Proof.** The detailed deduction of formulae (3.1)–(3.2) is given in [19]. □

**Lemma 3.2.** *For a prolate (along z) spheroid (figure 4b), the formula takes the form* [19] (*cf.* [11, 51, 59]):

$$V_z(x, y, 0) = 4\pi\gamma\rho \frac{\varepsilon}{e^3}\left[\ln\left(p + \sqrt{1 + p^2}\right) - \frac{p}{\sqrt{1 + p^2}}\right] z_0, \qquad (3.3)$$

*where*

$$\varepsilon = c/a > 1, \quad e = \sqrt{\varepsilon^2 - 1} > 0, \quad p = q/\sqrt{t}, \quad q = ea/r, \quad r = \sqrt{(x - x_0)^2 + (y - y_0)^2 + z_0^2},$$
$$t = \left\{1 - q^2 + \sqrt{(1 - q^2)^2 + 4q^2[(x - x_0)^2 + (y - y_0)^2]/r^2}\right\}/2. \qquad (3.4)$$

**Proof.** The detailed deduction of formulae (3.3)–(3.4) is given in [19]. □

In case of a sphere ($\varepsilon = c/a = 1$), the formula for $z$-intensity takes the form:

$$V_z(x, y, 0) = \frac{4}{3}\pi\gamma\rho a^3 \frac{z_0}{r^3}. \qquad (3.5)$$

In this paper, we use only formulae (3.1)–(3.5) for $z$-intensity $V_z(x, y, 0)$ of oblate and prolate spheroids and that of a sphere.

The use of spheroids for solving the inverse gravimetry problem makes this approach by adjoining to the Sretenskii approach. Recall that according to the Sretenskii class models, a body possesses the mean plane *P* if any straight line perpendicular to this plane intersects the body surface only at two points on different sides of the plane *P*. If a geological body possesses a mean plane, its gravity center is inside the body, and its density is constant (and is given!), then the inverse problem (determining the body shape from the potential) has a *unique solution*. Our approach extends the Sretenskii approach because (see section 5) we do not assume that the density ρ is given, but include it into a number of the sought parameters (along with parameters $a$, $\varepsilon$, $x_0$, $y_0$, $z_0$). Furthermore, the uniqueness of the solution is provided by introducing constraints on the parameters and a regularization[1] (see also below).

## 4. Initial approximations for deposit parameters

Consider $m \geq 1$ bodies. The isolines of the z-intensity $V_z(x, y, 0)$ produced by the bodies on Earth's surface are shown in figure 5. Let $V_z(x_i, y_i, 0)$ be measured at $N$ points $i = 1, \ldots, N$. The number of bodies $m$ and their coordinates $(x_0, y_0)_1, \ldots, (x_0, y_0)_m$ can be determined from the contour pattern (isolines).

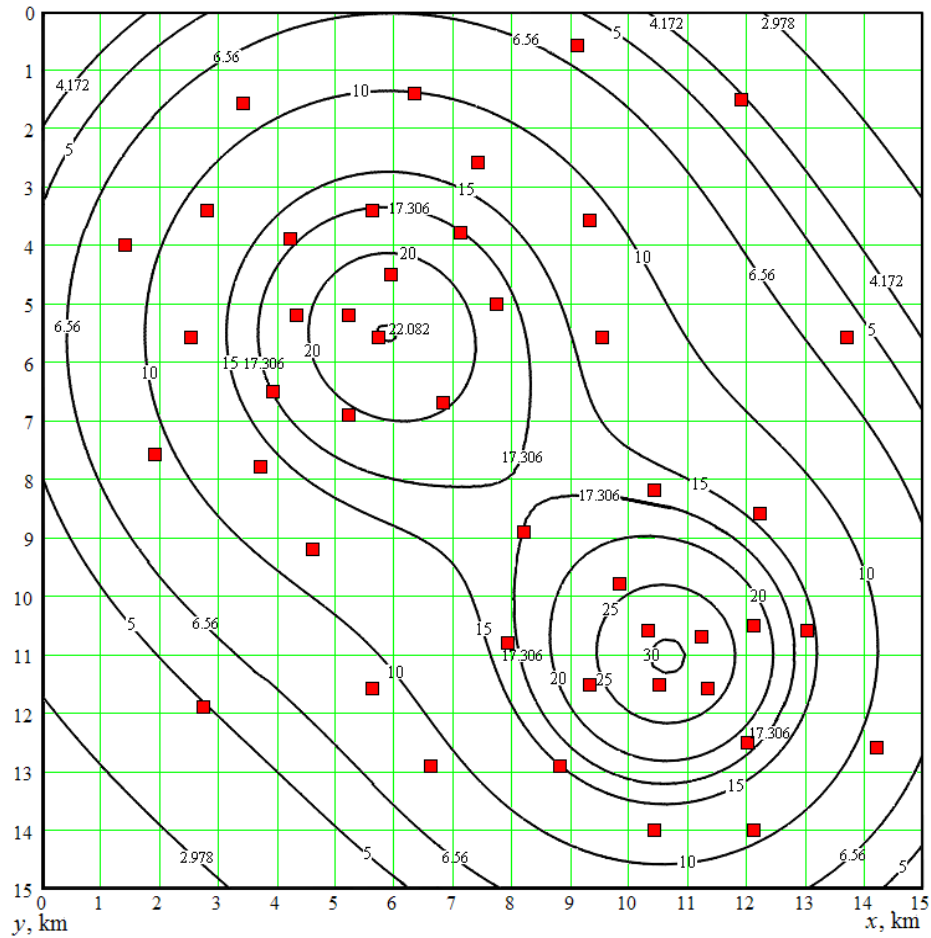

**Figure 5.** Isolines of intensity $V_z(x, y, 0)$, mGal from two bodies on frequent grids for *x* and *y*.

*4.1. Selection of bodies from isolines*

The following w a y of a selection of bodies from the contour pattern (isolines) is proposed. According to this way, *two conditions* must be fulfilled.

---

[1] If a solution is not unique, the Tikhonov regularization method chooses among a set of solutions a (unique) *normal pseudosolution* [1, 14, 28, 52, 56].

1. Valley in intensity $V_z$ between some two maxima (poles, peaks, hills) is not less than $\approx 20\%$ of intensities in poles (as in the Rayleigh criterion [45]).

2. The noise level $\delta V_z$ does not exceed 20% of the intensities in poles[2].

If both conditions are fulfilled, we assume that two peaks (and hence two bodies) are determined from isolines. For example, figure 5 shows two poles with intensities $V_{z1} \approx 22$ and $V_{z2} \approx 30$. We assume $V_z = (V_{z1}+V_{z2})/2 \approx 26$. The intensity between the poles is $v_z \approx 17$, i.e. the valley is $(V_z - v_z)/V_z = 0.346 \approx 35\%$. In this case, the noise level is $\delta V_z \approx 5\%$ (see section 6). As a result, both conditions are fulfilled and we can assume that two bodies are determined from isolines (see also figure 8 below).

The valley of 20% corresponds to some minimum distance between bodies in which they are delimited without mathematical processing. If the distance between the poles in figures 5 and 8 was $\approx 4$ km, the valley would be equal to $\approx 20\%$, i.e. 4 km is the limiting distance between the poles in which they are resolved. If the valley is less than $\approx 20\%$, the bodies can be resolved mainly mathematically. After separation of the bodies, we can make an estimate of some parameters of the deposit.

*4.2. Bulakh algorithm for estimating the depth and mass of a deposit*

In [6], an algorithm for estimating the depth $z_0$ and mass $M$ of each body is proposed in short form. We present this algorithm (the *Bulakh algorithm*) in more detail (cf. [19]). As a starting point, we assume that there is only one body, for example, the body matched by the $V_z$ contours in the upper left part of figure 5. We assume also that this body is a *homogeneous sphere*. Let us denote its center as $R(x_0, y_0, z_0)$, its mass as $M$, and the point at which $V_z = \max$ as $Q(x_0, y_0, 0)$ (figure 6). We suppose that $x_0$ and $y_0$ are known (for example, they are estimated by isolines). In Sections 5 and 6, values of $x_0$ and $y_0$ are refined.

The coordinates $x_0$ and $y_0$ are also the coordinates of the sphere center. Therefore, the unknowns will be only $z_0$ and $M$, whereas $V_{zP}(x,y,0)$ and $V_{zQ}(x_0,y_0,0)$ are measured.

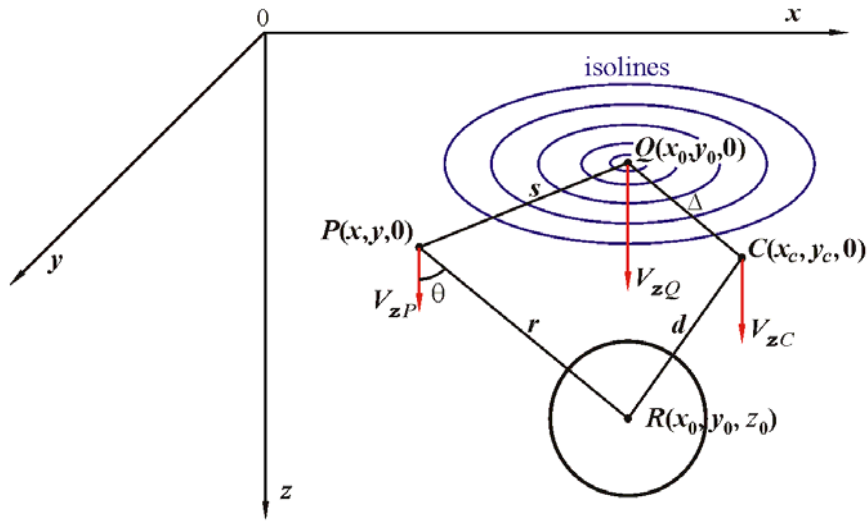

**Figure 6.** Deposit in the form of one body (homogeneous sphere).

---

[2]The value 20% may be changed. In addition, if the intensities in two poles are not equal: $V_{z1} \neq V_{z2}$, then we will use the value $V_z = (V_{z1}+V_{z2})/2$.

**Lemma 4.1.** *The z-intensity $V_z(x, y, 0)$ at point $P(x, y, 0)$ equals*

$$V_{zP} = \gamma \frac{M z_0}{(z_0^2 + s^2)^{3/2}}, \qquad (4.1)$$

*where $s = \sqrt{(x - x_0)^2 + (y - y_0)^2}$ is the distance between P and Q.*

**Proof.** We have

$$V_{zP} \equiv V_z(x, y, 0) = \gamma \frac{M}{r^2} \cos \theta = \gamma \frac{M}{r^2} \cdot \frac{z_0}{r},$$

where $r = \sqrt{z_0^2 + s^2}$ is the distance between *P* and *R*, and $\theta$ is the angle between $V_{zP}$ and *r*, which implies (4.1). □

**Corollary 4.1.** *The z-intensity at point Q equals*

$$V_{zQ} \equiv V_z(x_0, y_0, 0) = V_{z\max} = \gamma \frac{M}{z_0^2}. \qquad (4.2)$$

Relations (4.1) and (4.2) can be considered as a system of two equations with respect to $z_0$ and *M*. Using the notation $\nu = V_{zP}/V_{zQ}$, we obtain $\nu = (z_0/s)^3 / [(z_0/s)^2 + 1]^{3/2}$ or

$$\left( \frac{\mu^2}{\mu^2 + 1} \right)^{3/2} = \nu, \qquad (4.3)$$

where $\mu = z_0/s$.

Relation (4.3) is an equation with respect to $\mu = \mu(\nu)$ at given (measured) $\nu$. Its solution is

$$\mu(\nu) = \sqrt{\frac{\nu^{2/3}}{1 - \nu^{2/3}}}, \quad \nu \in (0, 1), \quad \mu \in (0, \infty). \qquad (4.4)$$

Dependence $\mu(\nu)$ is given in table 1.

**Table 1.** Dependence $\mu(\nu)$

| $\nu = \dfrac{V_z(s)}{V_{z\max}}$ | 0 | 0.1 | 0.2 | 0.3 | 0.4 | 0.5 | 0.6 | 0.7 | 0.8 | 0.9 | 1 |
|---|---|---|---|---|---|---|---|---|---|---|---|
| $\mu = \dfrac{z_0}{s}$ | 0 | 0.5240 | 0.7209 | 0.9011 | 1.0898 | 1.3048 | 1.5700 | 1.9301 | 2.4969 | 3.7071 | $\infty$ |

**Theorem 4.2** (The Bulakh algorithm). *Let a body be a homogeneous sphere. We assume that its coordinates $x_0$ and $y_0$ are known on a maximum intensity $V_{zQ}$ at the point $Q(x_0, y_0, 0)$ and a measurement of the intensity $V_{zP}$ is also performed at another point $P(x, y, 0)$. Then, the depth of the body equals $z_0 = \mu s$, where $\mu = \mu(\nu)$ is expressed by the formula (4.4), and $\nu = V_{zP}/V_{zQ}$, $s = \sqrt{(x - x_0)^2 + (y - y_0)^2}$. Furthermore, mass of the body equals*

$$M = \frac{1}{\gamma} \frac{(z_0^2 + s^2)^{3/2}}{z_0} V_{zP}. \qquad (4.5)$$

**Proof.** This follows from (4.1)–(4.4). □

**Remark 4.1.** It is preferable that the chosen point *P* is located away from the other bodies, e.g., at *x* < 5 km and *y* < 5 km in figure 5. One can make several estimates of μ (and $z_0$) from a few points *P* and average the result. This procedure should then be performed for each body.

**Remark 4.2.** The definition of depth $z_0$ and mass *M* of a body is possible only if $s \neq 0$, i.e. when measurements are performed at two different points *P* and *Q*. Moreover, to increase the accuracy of the algorithm, the following condition should be satisfied: $s \approx z_0$, which follows from the analysis of formulae (4.1)–(4.4) and numerical examples [6, 19].

*4.3. Generalization of the Bulakh algorithm*

The foregoing Bulakh algorithm can only be applied if the maximum intensity $V_{z\max} = V_z(x_0, y_0, 0)$ is measured at point $Q(x_0, y_0, 0)$ and the coordinates $x_0$ and $y_0$ are estimated (which is the case for the upper left body in figure 5). But if such a measurement is not made as in the lower right corner of figure 5, then we make a *generalization of the Bulakh algorithm*.

More concretely, we believe that a peak *Q* has been detected by isolines of intensity $V_z$ (lower right in figure 5) and its coordinates $x_0$ and $y_0$ are estimated, but the value of intensity $V_{zQ}$ in the peak is not directly measured. In addition to the points *P* and *Q*, we consider the point $C(x_c, y_c, 0)$ (figure 6), in which the value of $V_{zC}$ is measured, but not necessarily is equal to $V_{z\max}$. Furthermore, $V_{zC} > V_{zP}$, i.e. *C* is closer to the *Q*, than *P*.

The intensity $V_z$ at point $C(x_c, y_c, 0)$ equals (cf. (4.1))

$$V_{zC} = \gamma \frac{M z_0}{\left(z_0^2 + \Delta^2\right)^{3/2}}, \qquad (4.6)$$

where $\Delta = \sqrt{(x_c - x_0)^2 + (y_c - y_0)^2}$ is the distance between points *C* and *Q*.

The ratio of intensities at points *P* and *C* is (see (4.1) and (4.6))

$$\nu = \frac{V_{zP}}{V_{zC}} = \left(\frac{\mu^2 + \psi^2}{\mu^2 + 1}\right)^{3/2} \in (0,1),$$

where $\mu = z_0/s$, $\psi = \Delta/s$, or

$$\left(\frac{\mu^2(\nu, \psi) + \psi^2}{\mu^2(\nu, \psi) + 1}\right)^{3/2} = \nu. \qquad (4.7)$$

Relation (4.7) is an equation with respect to $\mu = \mu(\nu, \psi)$ at given ν and ψ. It generalizes equation (4.3). Its solution is (cf. (4.4))

$$\mu(\nu, \psi) = \sqrt{\frac{\nu^{2/3} - \psi^2}{1 - \nu^{2/3}}}, \qquad \nu \in (0,1), \quad \psi \in [0,1), \quad \mu \in (0, \infty). \qquad (4.8)$$

Dependence $\mu(\nu, \psi)$ is shown in figure 7.

After estimating μ by (4.8), we can estimate the *depth of the mass center of the body* (in the form of a sphere)[3]

$$z_0 = \mu s. \qquad (4.9)$$

Such estimate should be made for each body. Then we can estimate the body *mass M*.

---

[3] In Sections 5 and 6, the value $z_0$ is refined for a non-spherical (spheroidal) body.

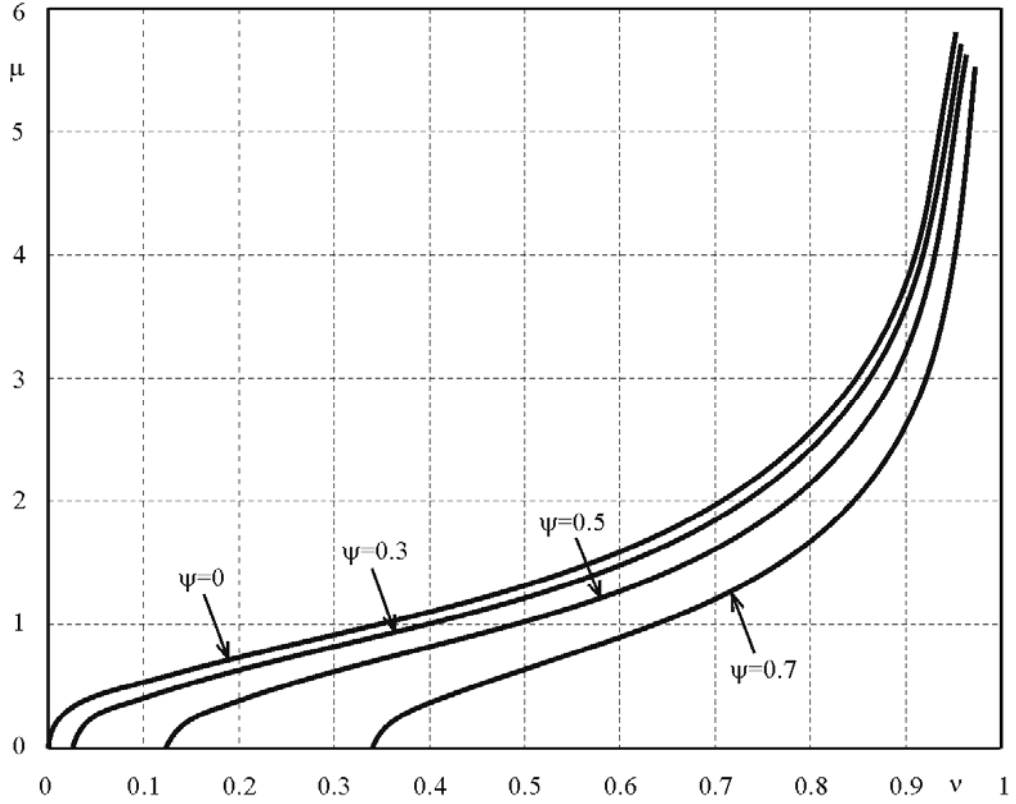

**Figure 7.** Dependence $\mu(\nu,\psi)$ (see (4.8)).

Furthermore, if $V_z$ is not measured at the point $(x_0, y_0, 0)$ of its maximum, but nevertheless, the coordinates $x_0, y_0$ are estimated, then one can use the measurement of $V_z$ at some point $C(x_c, y_c, 0)$ and obtain an estimate of the mass from (4.6):

$$M = \frac{1}{\gamma}\frac{(z_0^2 + \Delta^2)^{3/2}}{z_0}V_{zC}. \tag{4.10}$$

But if $V_z$ is measured at the point $(x_0, y_0, 0)$ of its maximum, one can obtain from (4.2)

$$M = \frac{1}{\gamma}z_0^2 V_{z\max}. \tag{4.11}$$

**Remark 4.3.** Formula (4.10) is similar to the formula (4.5) by writing. Nevertheless, these formulae differ essentially, namely, the point $P$ can not turn into $Q$ ($s$ is strictly greater than zero), and point $C$ may be at the maximum point $Q$ ($\Delta$ may be zero). Therefore, formulae (4.6)–(4.8) and (4.10) may be considered as a *generalization of the Bulakh algorithm*.

The obtained result can be formulated in the following theorem.

**Theorem 4.3** (generalization of the Bulakh algorithm). *Let a body be a homogeneous sphere. Let us assume that in some way (e.g., from isolines of intensity $V_z$), we select a local maximum (peak) $Q$ and estimate its coordinates $x_0$ and $y_0$, but does not measure the intensity $V_{zQ}$. Let us assume that we measure the intensity $V_{zP}$ at the point $P(x, y, 0)$ and $V_{zC}$ at the point $C(x_c, y_c, 0)$, more closed to $Q$. Then the body depth is $z_0 = \mu s$, where $\mu = \mu(\nu, \psi)$ is expressed by the formula (4.8), and $\nu = V_{zP}/V_{zC}$, $\psi = \Delta/s$, $s = \sqrt{(x - x_0)^2 + (y - y_0)^2}$, $\Delta = \sqrt{(x_c - x_0)^2 + (y_c - y_0)^2} \in [0, s)$, and the body mass $M$ is expressed by the formula (4.10). But if $\Delta = 0$ ($C$ turns into $Q$), then this theorem turns into theorem 4.2.*

**Remark 4.4.** Let the coordinates and distances $x$, $y$, $z$, $s$, $r$, $\Delta$, $d$ be expressed in kilometers, the mass $M$ in billions of tons, and the gravitational intensity (anomaly) $V_z$ in mGal. Then formulae (4.10) and (4.11) take the form

$$M = 0.15 \frac{(z_0^2 + \Delta^2)^{3/2}}{z_0} V_{zC}, \tag{4.12}$$

$$M = 0.15 z_0^2 V_{z\max}. \tag{4.13}$$

These results generalize the algorithm presented in [6] and yield a good initial approximation to $z_0$ and $M$ (as well as to $x_0$ and $y_0$) for each body. In the next section, we demonstrate how to make more precise these (and other) parameters of non-spherical (spheroidal) bodies.

## 5. Refinement of the deposit parameters

Consider the problem of *refining the parameters of spheroids* which approximate the bodies of a deposit. Let $\widetilde{V}_{zi} \equiv \widetilde{V}_z(x_i, y_i, 0)$, $i = 1, \ldots, N$ be the values of $V_z$ measured (with errors) at a series of points $(x_i, y_i, 0)$, $i = 1, \ldots, N$ on Earth's surface (in particular, on a ship), where $N$ is the number of measurement points. Let the values of $V_{zi} \equiv V_z(x_i, y_i, 0)$ are calculated using formulae (3.1), (3.3) or (3.5). Let $k$ be the number of sought parameters for each of the $m$ spheroids (the types of parameters are listed below), in total $mk$ parameters. Let us denote the deposit parameters as $p_j$, $j = 1, \ldots, mk$.

### 5.1. Application of the Tikhonov regularization method with constraints

The determination of the deposit parameters is an *ill-posed problem*, viz., unstable and having, generally speaking, a non-unique solution [3, 5, 6, 8, 17–19, 26, 30, 33, 46–49]. Even if a body is star-shaped and has the Sretenskii mean plane, the solution is anyway non-unique if the body density ρ is not given. A quite effective way to eliminate the non-uniqueness of the solution is to use *constraints* on the solution with the Tikhonov regularization method.

We will solve this problem by *minimizing the Tikhonov smoothing* (*parametric*) *functional* [56] (cf. [3, 10, 15, 16, 24, 25, 30, 55, 61]) in two variants [17–19]:

$$F_1 \equiv \sum_{i=1}^{N} \left(\widetilde{V}_{zi} - V_{zi}\right)^2 + \alpha \sum_{j=1}^{mk} q_j (p_j - p_{\mathrm{mid}\,j})^2 = \min_{p_1, \ldots, p_{mk}} \tag{5.1}$$

or

$$F_2 \equiv \sum_{i=1}^{N} \left(\widetilde{V}_{zi} - V_{zi}\right)^2 + \alpha \sum_{j=1}^{mk} q_j p_j^2 = \min_{p_1, \ldots, p_{mk}}, \tag{5.2}$$

where $q_j$ are the weights and $\alpha > 0$ is the regularization parameter (the Lagrange multiplier).

Functional (5.1) or (5.2) contains the discrepancy (or the data misfit) $\sum_{i=1}^{N}(\widetilde{V}_{zi} - V_{zi})^2$ and stabilizer $\sum_{j=1}^{mk} q_j (p_j - p_{\mathrm{mid}\,j})^2$ or $\sum_{j=1}^{mk} q_j p_j^2$. In works [40, 60] et al., different variants of the stabilizer are considered conformably to the geophysical inversion problem with the purpose of focusing the noisy and insufficiently clear geophysical inverse images [10, 55, 61]. In [29], the problem of minimizing the objective function using discrepancy and weight functions is considered. Furthermore, the objective function includes a priori information and additional geophysical and geological constraints. The minimization of the objective function is achieved by using a generalized subspace inversion algorithm.

In works [15, 22] et al., the inverse gravimetry problem is solved by a sparse regularization oriented on solving under-determined or ill-conditioned systems of linear equations and based on a matching pursuit algorithm.

In order to enhance the stability and eliminate the solution non-uniqueness, we impose *constraints on the solution* by setting the ranges of values for the parameters (cf. [16, 49, 62]):

$$p_{\min j} \leq p_j \leq p_{\max j}, \quad j = 1, \ldots, mk, \tag{5.3}$$

where $p_{\min j}$ and $p_{\max j}$ must be estimated from an additional (a priori) information[4].

We assume for the weights that $q_j = 1/p_{\mathrm{mid}\,j}^2$, where $p_{\mathrm{mid}\,j} = (p_{\min j} + p_{\max j})/2$. Such a problem is a problem of nonlinear programming with regularization [21, 35, 54, 56].

**Remark.** We introduce the weights $q_j$ because the sought parameters $p_j$ can have a different physical dimension and a different order of magnitudes, and the introduction of the weights makes the summands $q_j(p_j - p_{\mathrm{mid}\,j})^2$ in (5.1) or $q_j p_j^2$ in (5.2) dimensionless quantities which are close to each other.

The mean square error of the solution equals [17, 19]

$$\delta = \left[ \frac{1}{mk} \sum_{j=1}^{mk} q_j (p_j - \bar{p}_j)^2 \right]^{1/2}, \tag{5.4}$$

where $p_j = p_j(\alpha)$, $\delta = \delta(\alpha)$, and $\bar{p}_j$ are 'exact' $p_j$ estimated from the modelled body.

## 5.2. Scheme for calculating the parameters of deposit

We suggest the following s c h e m e  for calculating the parameters of deposit bodies.

1. The number of bodies $m$ and the values of $x_0$, $y_0$ for each body are estimated from the $V_z$ isolines pattern (see figure 5 or figure 10 further).

2. The average value of $\mu$ is determined using formula (4.4) or/and (4.8) and $z_0 = \mu s$ is estimated for each body from the measurements of $V_z$ at the series of points.

3. The mass $M$ of each body is estimated using formulae (4.5), (4.10)–(4.13).

4. The $k = 5$ parameters $p_1 = \varepsilon$, $p_2 = \rho$, $p_3 = x_0$, $p_4 = y_0$, and $p_5 = z_0$ are sought for each body, while estimations of $x_0, y_0, z_0$ obtained by the Bulakh algorithm are used as the initial approximations assuming $p_{\mathrm{mid}\,3} = x_0$, $p_{\mathrm{mid}\,4} = y_0$, $p_{\mathrm{mid}\,5} = z_0$. Furthermore, the estimate of $M$ is used in the calculation of the semiaxis $a$,

$$a = \sqrt[3]{\frac{M}{(4/3)\pi\varepsilon\rho}}. \tag{5.5}$$

## 5.3. A way of decremental constraints

A number of methods is developed for the problem of minimizing functionals with constraints (the methods of conditional gradient, conjugate gradients, steepest descent, ravines, chords, et al.) [21, 54, 56]. In exploration geophysics, the *method of gradient descent* is often used (the works of Kantorovich, et al.) [5]. However, the point of convergence in this method strongly depends on the initial approximation of the solution. Besides that, this method is using the derivatives of a functional and is not using constraints on solution.

---

[4] In Section 6.1, it is indicated that the constraints can be estimated by the results of preliminary study of the region. We note also that in [23], a priori information is formulated in the form of a fuzzy set in the parameter space. In such a formulation the inverse problem becomes a multiobjective optimization problem with two objective functions, one of them is a membership function of the fuzzy set of feasible solutions, the other is the conditional probability density function of the observed data.

In this paper, for minimizing the functional (5.1) or (5.2), we use a modification of the *coordinate descent method* [13] which is, to a certain degree, effective in case when constraints have form (5.3). A special feature of this method is that it does not allow the solution to fall outside the constraints (5.3), thus providing stability and convergence of the solution within the interval ('corridor') given by inequalities (5.3). However, it is not always possible to set a sufficiently narrow interval. In this case, first, one should set wide constraints, and then one must do more narrow constraints iteratively so that the solution does not go out of the interval (5.3). The following examples are solved in exactly this way. We call this method a *way of decremental constraints*.

## 6. Numerical examples (modelling of the direct and inverse problems)

*Program packages* IGP2 and IGP5 (Inverse Gravimetry Problem, var. 2 and 5) were developed for modelling the direct and inverse gravimetry problems for the case when a deposit consists of two or five bodies ($m = 2$ or $m = 5$). The calculations were implemented in MS Fortran PowerStation 4.0 (Fortran90) and graphs were created in MathCAD, CorelDRAW, and Paint. The two following *numerical examples* were solved.

### 6.1. Example 1

The deposit was given in the form of two bodies (figures 1 and 2).

**The direct problem.** In the *direct problem* (modelling of the intensities $V_z$ created by the deposit), there were $N = 45$ points $(x_i, y_i, 0)$, $i = 1, \ldots, N$ for measuring $V_z$ on Earth's surface (the small squares in figure 5). The intensity $V_z(x_i, y_i, 0)$ was simulated at each point $(x_i, y_i, 0)$, $i = 1, \ldots, N$, by summing the elementary bars for the two bodies using (2.1) with the discretization step $dx' = dy' = 0.25$ km for the first body and $dx' = dy' = 0.125$ km for the second one having $z'_{\min}(x', y')$ and $z'_{\max}(x', y')$ given (figure 3).

Normal errors $\delta V_z$ with zero expectation and standard deviation of $\sigma = 1$ mGal (which is $\approx 6\%$ of the average value of $V_z$ and $\approx 4\%$ of $V_{z\,\max}$) were generated by a random-number generator [44] and added to the exact values $V_z(x_i, y_i, 0)$. The errors $\delta V_z$ model measurement errors, small-scale non-homogeneities of a medium, etc.

$V_z(x_i, y_i, 0)$ were plotted as isolines using command Contour Plot in MathCAD. In figure 5, the isolines of $V_z(x, y, 0)$ (mGal) are given at $x \in [0,15]$ km, $y \in [0,15]$ km. The density of the first body was $\rho = 1.6\,\text{g}/\text{cm}^3$ and that of the second one was $\rho = 2.6\,\text{g}/\text{cm}^3$. Such density anomalies can be associated, e.g., with the upper ore body and shales, respectively.

However, in practice, continuous data on $V_z(x, y, 0)$ are usually unknown, whereas (noisy) values of $\widetilde{V}_z$ are known only at measurement points $(x_i, y_i, 0)$, $i = 1, \ldots, N$. In figure 8, the isolines of $z$-intensity $\widetilde{V}_z = V_z + \delta V_z$ (mGal) are plotted from the measured values of $\widetilde{V}_z(x_i, y_i, 0)$, $i = 1, \ldots, N$.

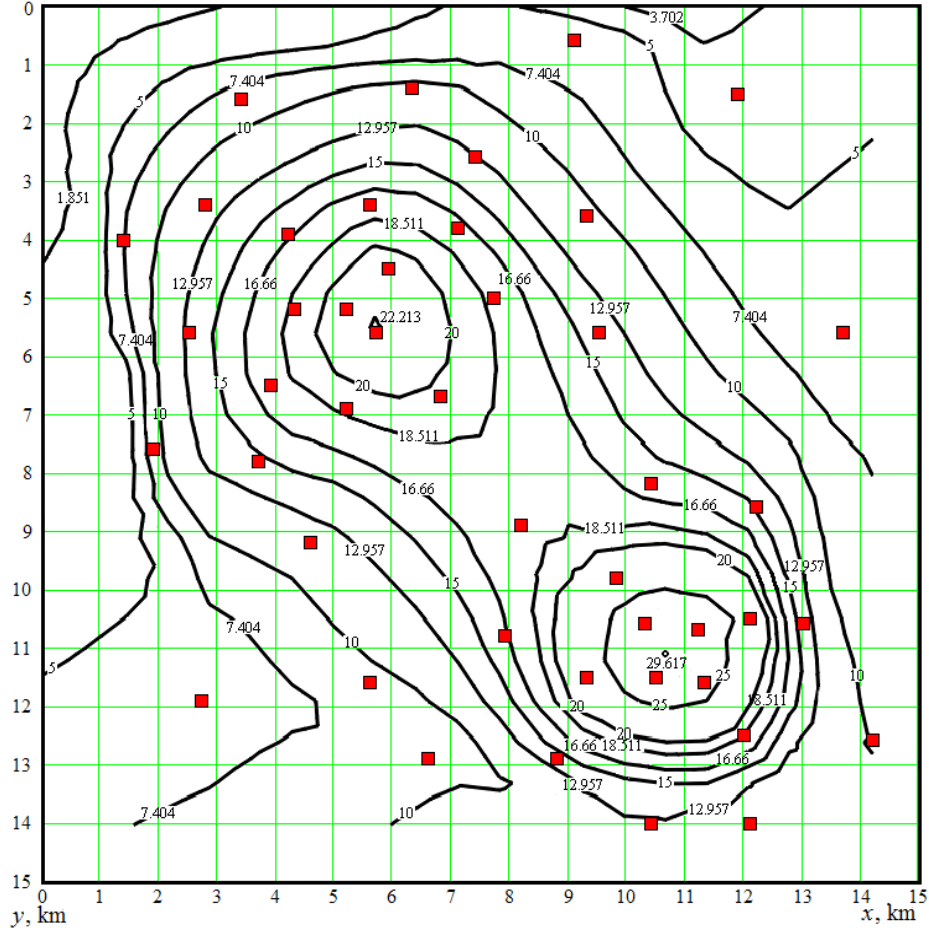

**Figure 8.** The isolines of the intensity $\widetilde{V}_z(x, y, 0)$, mGal, induced by two bodies. The isolines are plotted using $N = 45$ measurement points (example 1).

**Remark.** Strictly speaking, isolines of $V_z(x, y, 0)$ are not necessary for constructing deposit model. As we see from (5.1) and (5.2), the discrete values of $\widetilde{V}_{zi}$ are only used. Nevertheless, isolines are necessary for estimating the number of bodies and their coordinates $x_0$ and $y_0$.

**The inverse problem.** We present in the following the results of solving the *inverse problem*, viz. the determination of the parameters of the two deposit bodies.

First, $x_0$, $y_0$ were evaluated for both bodies from the isolines in figure 8. For the first body, $x_0 = 5.7$ and $y_0 = 5.6$ km were obtained. The estimates of $z_0 = \mu s$, which were obtained according to formulae (4.4) and (4.8) for few points $P$ and $Q$ located in the area $x \leq 7$ & $y \leq 6$ km, vary from 3 to 5.5 km, or on average $z_0 = 4.9$ km. The estimates of mass $M$ obtained by formulae (4.10) and (4.11) (or (4.12) and (4.13)) vary from 30 to 105 bln t, or on average $M = 67.5$ bln t.

For the second body, $x_0 = 10.65$ and $y_0 = 11.1$ km. The estimates of $z_0 = \mu s$ according to formula (4.8) for few points at $x \geq 8$ & $y \geq 10$ km vary from 2.9 to 3.7 km, or on average $z_0 = 3.3$ km. The estimates of mass $M$ obtained by formula (4.13) vary from 37 to 60 bln t, or on average $M = 48.5$ bln t.

Then, the following ten sought parameters were refined by minimizing the functional $F_1$ or $F_2$ (see (5.1) or (5.2)): $p_1 = \varepsilon$, $p_2 = \rho$, $p_3 = x_0$, $p_4 = y_0$, $p_5 = z_0$ for the first body and $p_6 = \varepsilon$, $p_7 = \rho$, $p_8 = x_0$, $p_9 = y_0$, $p_{10} = z_0$ for the second one.

In table 2, the initial lower constraints $p_{\min}$, the upper constraints $p_{\max}$ and the mean values $p_{\mathrm{mid}} = (p_{\min} + p_{\max})/2$ are presented. These (approximate) values of the constraints can be estimated by the experts in this field on the basis of additional data obtained by results of preliminary investigations of the region[5]. Constraints on $\rho$ are especially important, since the density knowledge connects with the problem of uniqueness. Then these values are determined more accurately such that the sought solution $p_i, \ldots, p_{10}$ does not fall outside the limits of the interval (5.3). Furthermore, estimate of the semiaxis $a$ is realized by formula (5.5).

**Table 2.** Refined values of the constraints on the parameters and the obtained solution

| $p$ | $\varepsilon$ | $\rho$, g/cm$^3$ | $x_0$, km | $y_0$, km | $z_0$, km | $v$, km$^3$ | $M$, bln t | $a$, km |
|---|---|---|---|---|---|---|---|---|
| The first body ||||||||||
| $p_{\min}$ | 0.2 | 1.1 | 5.4 | 5.2 | 4.0 | | | |
| $p_{\max}$ | 0.6 | 1.7 | 6.0 | 6.0 | 5.8 | | | |
| $p_{\mathrm{mid}}$ | 0.4 | 1.4 | 5.7 | 5.6 | 4.9 | | | |
| solution | 0.595 | 1.69 | 5.75 | 5.48 | 4.40 | 39.86 | 67.50 | 2.52 |
| exact values | 0.51 | 1.6 | 5.7 | 5.3 | 4.2 | 39.67 | 63.48 | 2.75 |
| The second body ||||||||||
| $p_{\min}$ | 1.8 | 2.3 | 10.3 | 10.2 | 2.3 | | | |
| $p_{\max}$ | 2.2 | 2.9 | 11.0 | 12.0 | 4.3 | | | |
| $p_{\mathrm{mid}}$ | 2.0 | 2.6 | 10.65 | 11.1 | 3.3 | | | |
| solution | 2.035 | 2.65 | 10.71 | 11.90 | 3.81 | 18.28 | 48.50 | 1.29 |
| exact values | 1.96 | 2.6 | 10.7 | 11.1 | 3.8 | 19.06 | 49.55 | 1.375 |

The dependence of the mean square error $\delta$ of the solution (see (5.4)) on the regularization parameter $\alpha$ was calculated. Furthermore, values of $\varepsilon$, $\rho$, $x_0$, $y_0$, $z_0$, $v$, $M$ and $a$ based on the deposit estimation (figures 1 and 2) and given in Table 2, were taken as 'exact' parameters $\bar{\varepsilon}$, $\bar{\rho}$, $\bar{x}_0$, $\bar{y}_0$, $\bar{z}_0$, $\bar{v}$, $\bar{M}$ and $\bar{a}$ of the deposit bodies.

Note that the exact volumes $\bar{v}$ and masses $\bar{M}$ of bodies were calculated by summation over the bars, and the values of $\bar{\varepsilon}$, $\bar{\rho}$, $\bar{x}_0$, $\bar{y}_0$ and $\bar{z}_0$ were estimated by the form of bodies. For example, the $x$ half-length of the first body is 2.5 km and the $y$ half-length is 3 km, i.e. the average half-length in the plane $xy$ is $\bar{a} = 2.75$ km; the $z$ half-length is 1.4 km, i.e. the ratio between the $z$ and $xy$ half-lengthes is $\bar{\varepsilon} = 0.51$, etc.

Since a regularization is used, it is necessary to draw our attention to the problem of choosing the value of the regularization parameter $\alpha$. A number of approaches for choosing $\alpha$ is developed (the Morozov discrepancy principle, L-curve algorithm, generalized cross-validation method et al.) [1, 14, 15, 28, 44, 55, 56]. However, in example 1, the rigid con-straints are imposed on the solution $p$ and they provide stability and uniqueness of the solution without regularization. Practically, the problem has the same solution for any $\alpha \gtrsim 10^{-8}$ and the solution does not depend on the type of the functional ($F_1$ or $F_2$). But if the constraints are less rigid than those given in table 2, the error $\delta(\alpha)$ and the solution $p(\alpha)$ depend on $\alpha$, and the value of the error $\delta$ at some $\alpha = \alpha_{opt}$ can be smaller than that at $\alpha = 0$.

---

[5] As a result of such investigations, one can determine which bodies (ore, oil, quartzites, shales, etc.) occur in the region. This makes possible to estimate their densities. As to $\varepsilon$, then at first we can put $p_{\mathrm{mid1}} = \varepsilon = 1$.

In figure 9, the pictures of the two calculated spheroids are given. It can be seen that the modelling of a deposit by spheroids yields a quite satisfactory result.

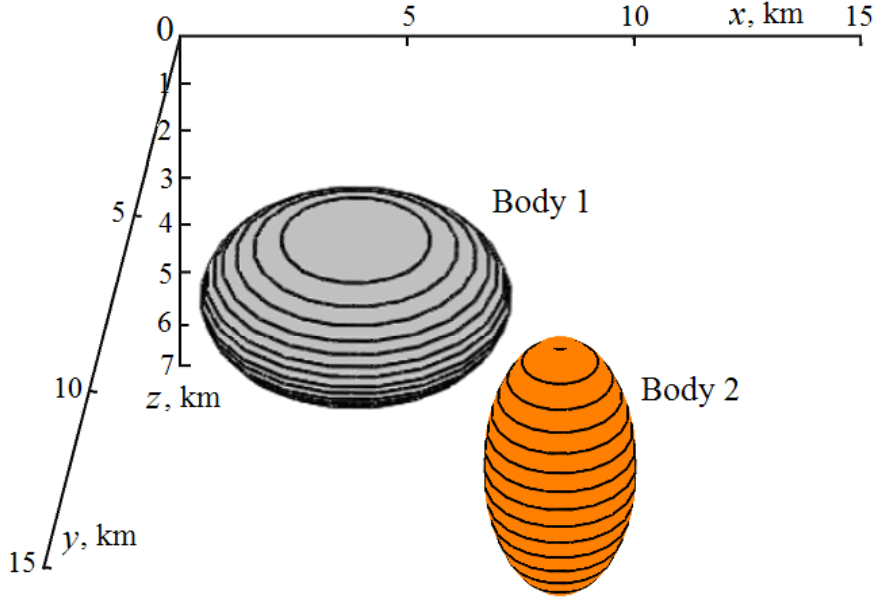

**Figure 9**. Two spheroids imitating the deposit in example 1.

It should be pointed out that the initial value of the functional $F_1$ (before its minimization) equals $F_1 = 0.167 \cdot 10^{-2}$, and after the minimization by the coordinate descent method $F_1 = 0.197 \cdot 10^{-3}$, i.e. the value of $F_1$ decreased by one order of magnitude.

### 6.2. Example 2

In the capacity of the second example, we could consider an example based on real measured fields on Earth's surface. However, in this case the exact values of the parameters of the bodies would not be known and it would be difficult to estimate the errors in determining the parameters. Of course, instead of estimating and comparing the parameters, one can use a comparison of the measured isolines of $\tilde{V}_z$ (type of figures 8 and 10) and the isolines obtained as a result of calculating the spheroids (see figure 11 later). However, it will be an indirect estimating the errors of body parameters. In addition, because of ill-posedness of the problem a small difference of isolines does not guarantee a small difference of body parameters.

It could be considered a *synthetic example* where the real data about deposit bodies are used, then calculations of the fields $\tilde{V}_z(x,y,0)$ are performed, after which the inverse problem of determining the parameters of the spheroids is solved. In this case, we will be able to compare the calculated parameters of the spheroids with the parameters of the deposit bodies. However, it is difficult to get real data about deposit bodies with the necessary accuracy.

Therefore, as the second example we have considered an example in which (as in example 1) the parameters of the deposit bodies are given (see table 3 below), but the deposit has more bodies than in example 1. Further, the intensities $V_z$ are calculated at $N = 73$ points on Earth's surface (squares in figure 10) and 3% errors are added to values of $V_z$. Example 2 is more complicated than example 1, and its purpose is to show potential possibilities of the method on increasing the number of bodies.

Based on the values of $\widetilde{V}_z$ using computer graphics of the system MathCAD, the $z$-intensity contours were constructed (figure 10).

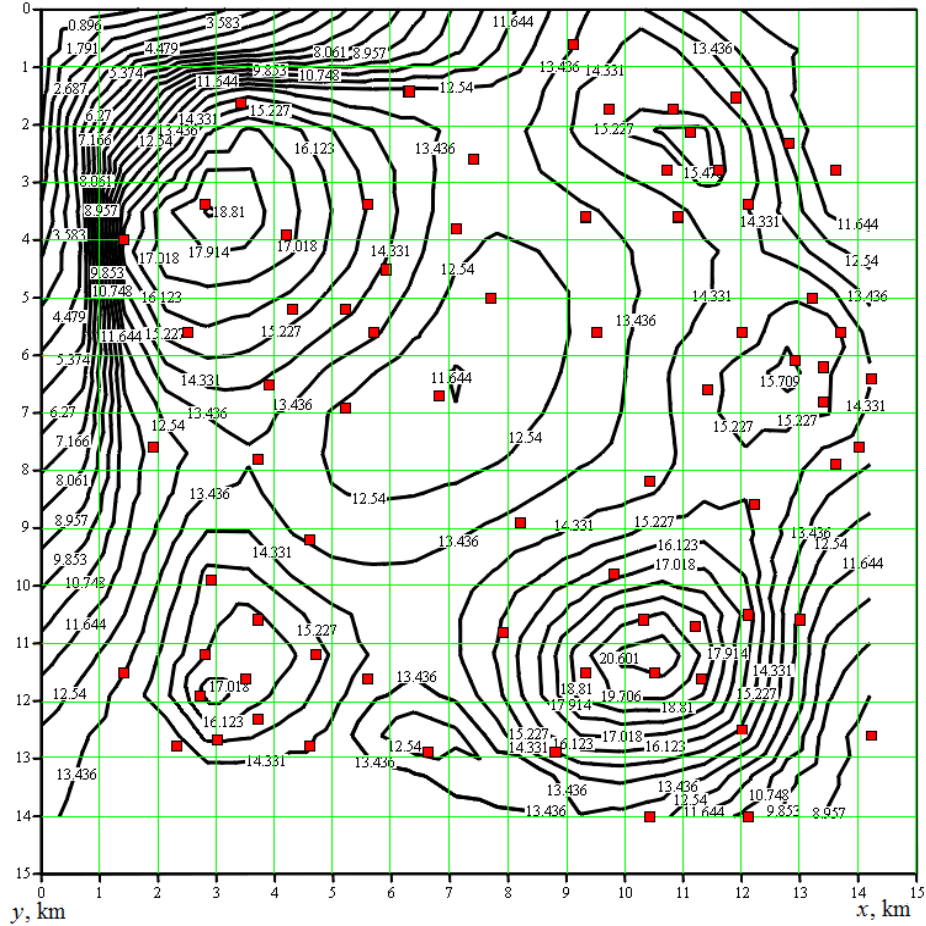

**Figure 10.** The isolines of the intensity $\widetilde{V}_z(x,y,0)$, mGal, induced by more than two bodies. The isolines are plotted with the help of computerized graphics of the system MathCAD using $N = 73$ measurement points (example 2).

By the isolines of figure 10 using the way of a selection of bodies from the isolines (see section 4.1), we conclude that the number of bodies in deposit is $m = 5$. Furthermore, it is also possible to estimate the coordinates $(x_0, y_0)$ for each body. Using the generalization of the Bulakh algorithm and formulae (4.4), (4.8), (4.12) and (4.13), the estimations of the depths $z_0 = \mu s$ and masses $M$ of the bodies were obtained as well. The results are as follows ($x_0, y_0, z_0$ in km, $M$ in bln t):

body 1: $x_0 = 3$, $y_0 = 3.5$, $z_0 = 4.45$, $M = 55.9$;

body 2: $x_0 = 10.3$, $y_0 = 11.3$, $z_0 = 4.1$, $M = 47.5$;

body 3: $x_0 = 3$, $y_0 = 11.8$, $z_0 = 3.8$, $M = 37.5$;

body 4: $x_0 = 11$, $y_0 = 2$, $z_0 = 4.4$, $M = 43.1$;

body 5: $x_0 = 12.9$, $y_0 = 6.3$, $z_0 = 5.1$, $M = 30$.

Then, by minimizing the functional $F_1$ (see (5.1)) the sought parameters of the bodies-spheroids $p_1 = \varepsilon$, $p_2 = \rho$, $p_3 = x_0$, $p_4 = y_0$, $p_5 = z_0$ for body 1 and similarly for bodies 2–5 (in total, $mk = 25$ parameters) were determined more accurately. And at that, first, more wide constraints on the values of the parameters were used, and then, they were narrowed step by step such that the sought solution $p_1,\ldots,p_{25}$ did not fall out of the interval (5.3). After refinement of $z_0$, $M$ is refined according to (4.13).

The regularization parameter was chosen to be $\alpha = 10^{-8}$. However, its value has very weak influence on the solution due to the constraints on the solution already act as regularization. The obtained solution is summarized in table 3.

**Table 3**. Obtained ($\widetilde{s}$) and exact ($\overline{s}$) solutions for example 2

| Body | | $a$, km | $\varepsilon$ | $\rho$, g/cm$^3$ | $x_0$, km | $y_0$, km | $z_0$, km | $\upsilon$, km$^3$ | $M$, bln t |
|---|---|---|---|---|---|---|---|---|---|
| 1 | $\widetilde{s}$ | 2.341 | 0.54 | 1.64 | 2.84 | 3.34 | 4.40 | 29.02 | 47.60 |
|   | $\overline{s}$ | 2.5   | 0.51 | 1.6  | 2.8  | 3.3  | 4.2  | 29.81 | 47.69 |
| 2 | $\widetilde{s}$ | 1.306 | 1.60 | 2.34 | 10.38 | 12.00 | 4.00 | 14.91 | 34.90 |
|   | $\overline{s}$ | 1.375 | 1.56 | 2.3  | 10.3  | 11.7  | 3.8  | 15.17 | 34.89 |
| 3 | $\widetilde{s}$ | 1.696 | 1.04 | 1.54 | 2.84  | 12.00 | 4.20 | 21.23 | 32.70 |
|   | $\overline{s}$ | 1.8   | 1.0  | 1.5  | 2.8   | 11.8  | 4.0  | 21.81 | 32.72 |
| 4 | $\widetilde{s}$ | 1.235 | 1.44 | 2.74 | 11.00 | 1.24  | 4.60 | 11.35 | 31.10 |
|   | $\overline{s}$ | 1.3   | 1.4  | 2.7  | 10.8  | 1.2   | 4.4  | 11.51 | 31.06 |
| 5 | $\widetilde{s}$ | 1.320 | 0.74 | 3.34 | 13.70 | 6.50  | 4.10 | 7.13  | 23.80 |
|   | $\overline{s}$ | 1.4   | 0.7  | 3.3  | 13.5  | 6.3   | 3.9  | 7.18  | 23.71 |

Based on the values of the parameters of the five bodies-spheroids summarized in table 3, the $z$-intensities $V_z(x, y, 0)$, $x \in [0,15]$ km, $y \in [0,15]$ km were calculated and the isolines $V_z$ were plotted (figure 11).

The comparison between measured $\widetilde{V}_z$ (figure 10) and calculated $V_z$ (figure 11) isolines shows that the main goal in this numerical example was reached despite it was complicated enough for processing, namely, the five bodies were identified; their surface coordinates $x_0, y_0$ were estimated; the values of $V_z$ on the five calculated 'poles' (figure 11) correspond to the values of $\widetilde{V}_z$ on the measured 'poles' (figure 10), etc.

In figure 12, the images of the five calculated spheroids are given. They can be identified with a deposit consisting of several bodies (ore, oil, intrusion, quartzites, shales, granite, basalt, etc.).

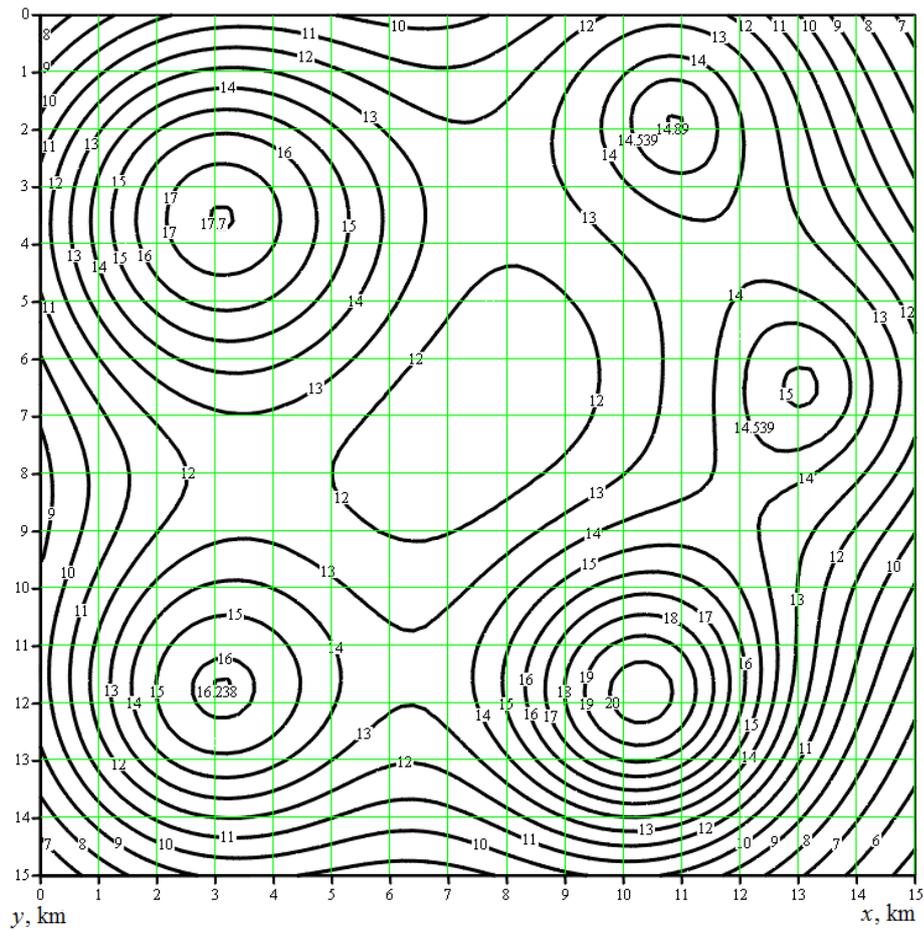

**Figure 11.** The isolines of the intensity $V_z(x, y, 0)$, mGal, induced by five spheroids (example 2)

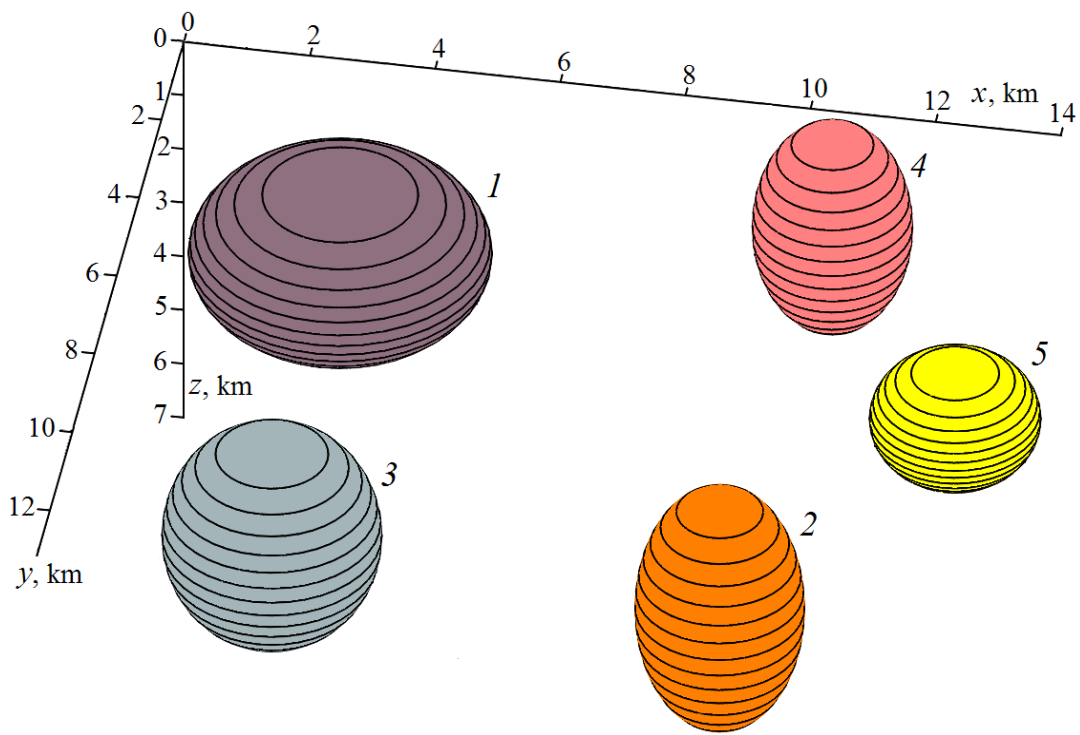

**Figure 12.** The five spheroids imitating the deposit in example 2.

## 7. Conclusion

Based on the results obtained in the present paper as well as in [17–19], the following c o n -
c l u s i o n s  were made:

1. Constraints on the solution (the parameters of the deposit model, e.g., in the form of inequalities (5.3)) are helpful in eliminating the non-uniqueness of the inverse gravimetry problem and make it stable. The Tikhonov regularization method is efficient in this case as well (see (5.1), (5.2)).

2. It is possible to obtain a stable (and unique) solution to the inverse gravimetry problem without a regularization, but by imposing rigid constraints on the solution.

3. The computerized implementation of the functional minimization problem for $F_1$ and $F_2$ (see (5.1) or (5.2)) requires solely one-dimensional arrays ($V_{zi}$ and others), but not matrices (cf. [49]). The time required for solving the functional minimization problem by the coordinate descent method in the above given examples was of the order of 1–10 s at a CPU clock rate of about 1 GHz. Hence, the proposed technique does not require large computer resources (memory and time).

4. T h e  m a i n  c o n c l u s i o n. The use of spheroids in deposit modelling as well as the application of the smoothing functional minimization method with constraints on the sought parameters yields to solve the inverse gravimetry problem enough accurately, steady, uniquely, visually and with little expenditure of computer memory and time even in case of significant errors in the initial data.

## Acknowledgements

This work is supported by the Russian Foundation for Basic Research (grants no. 09-08-00034 and 13-08-00442); Energinet.dk (ForskEL 2009-1-10246, ForskEL 2008-1-0079); the European Social Fund (ESF) in the Czech Republic, the Ministry of Education, Youth and Sports (MEYS) of the Czech Republic (POST-UP CZ.1.07/2.3.00/30.0004); the institutional fund of Palacký University Olomouc.